\documentclass[12pt,a4paper]{amsart}
\usepackage{a4wide}
\usepackage[english,russian]{babel}
\usepackage[T2A]{fontenc}
\usepackage[cp1251]{inputenc} 
\usepackage{amsfonts}
\usepackage{amssymb, amsthm, amscd}
\usepackage{amsmath}
\usepackage{mathtools}
\usepackage{needspace}
\usepackage{etoolbox}
\usepackage{lipsum}
\usepackage{comment}
\usepackage{cmap}
\usepackage[pdftex]{graphicx}
\usepackage[unicode]{hyperref}
\usepackage[matrix,arrow,curve]{xy}
\usepackage[usenames,dvipsnames]{xcolor}
\usepackage{colortbl}
\usepackage{textcomp}
\usepackage{cite}
\usepackage{euscript}

\newcommand{\N}{\mathbb{N}}

\newcommand{\F}{\mathbb{F}}

\newcommand{\eps}{\varepsilon}

\renewcommand{\ge}{\geqslant}
\renewcommand{\le}{\leqslant}
\newcommand{\sm}{\setminus}
\newcommand{\map}[3]{#1\colon #2\to #3}

\theoremstyle{plain}
\newtheorem{thm}{Theorem}
\newtheorem{lem}{Lemma}

\newtheorem{cor}{Corollary}

\theoremstyle{definition}

\theoremstyle{remark}

\AtBeginEnvironment{thm}{\begin{samepage}}
\AtEndEnvironment{thm}{\end{samepage}}
\AtBeginEnvironment{lem}{\begin{samepage}}
\AtEndEnvironment{lem}{\end{samepage}}
\AtBeginEnvironment{st}{\begin{samepage}}
\AtEndEnvironment{st}{\end{samepage}}
\AtBeginEnvironment{crit}{\begin{samepage}}
\AtEndEnvironment{crit}{\end{samepage}}
\AtBeginEnvironment{ax}{\begin{samepage}}
\AtEndEnvironment{ax}{\end{samepage}}
\AtBeginEnvironment{defn}{\begin{samepage}}
\AtEndEnvironment{defn}{\end{samepage}}
\AtBeginEnvironment{cor}{\begin{samepage}}
\AtEndEnvironment{cor}{\end{samepage}}
\AtBeginEnvironment{note}{\begin{samepage}}
\AtEndEnvironment{note}{\end{samepage}}
\AtBeginEnvironment{prop}{\begin{samepage}}
\AtEndEnvironment{prop}{\end{samepage}}

\newcommand{\GL}{\mathop{\mathrm{GL}}\nolimits}

\newcommand{\tc}{\text{,}}
\newcommand{\tp}{\text{.}}

\renewcommand{\tilde}{\widetilde}

\DeclareMathOperator{\diag}{diag}
\DeclarePairedDelimiter{\ceil}{\lceil}{\rceil}

\makeatletter
\def\@settitle{\begin{center}%
    \baselineskip14\p@\relax
    \bfseries
    \@title
  \end{center}%
}

\def\@evenhead{\hfil\sc p. gvozdevsky\hfil}
\def\@oddhead{\hfil\sc commutator lengths \hfil}
\makeatother

\title{Commutator lengths in general linear group over a skew-field.}

\address{Chebyshev Laboratory, St. Petersburg State University, 14th Line V.O., 29, Saint Petersburg 199178 Russia.}

\thanks{Research is supported by «Native towns», a social investment program of PJSC ''Gazprom Neft'', and also by grant given as subsidies from Russian federal budget for creation and development of international mathematical centres, agreement between MES and PDMI RSA  \textnumero 075-15-2019-1620 from November 8, 2019.}

\author{P. Gvozdevsky}

\begin{document}
\selectlanguage{english}

\maketitle

\begin{abstract}
	We give an upper and lower estimate for the maximal commutator length of a noncentral element of the elementary subgroup of the general linear group over a skew-field based on the maximal commutator length of an element of the multiplicative group of that skew-field. 
\end{abstract}

\section{Introduction}

Questions on commutator lengths are questions on specific word maps. In general case, for a group $G$ {\it a word map} is a map
$$
G^n\to G\tc\qquad
(g_1,\ldots,g_m)\mapsto w(g_1,\ldots,g_m)\tc
$$     
where $w$ is a word of $m$ variables (i.e an element of the free group on $m$ generators), and where $w(g_1,\ldots,g_n)$ is what we get if we substitute the elements $g_1$,$\ldots$,$g_n$ into this word.

By definition, the {\it commutator length} of the element $g\in[G,G]$ if not greater than $d$ if and only if $g$ belongs to the image of the word map given by the word $w=[x_1,y_1]\ldots[x_d,y_d]$ of $2d$ variables. In other words, if it can be expressed as a product of $d$ commutators. 

In recent studies of word maps, the case of $G=\mathcal{G}(K)$, where $\mathcal{G}$ is a simple or semisimple algebraic group defined over the field $K$, is of particular interest, see references in \cite{GordKunPlotRAN}, \cite{GordKunPlotJAlg},\cite{GordKunPlotJIJAC}. Almost all the results in this direction have been proved for the cases where $K$ is an algebraically closed field, or where $\mathcal{G}$ is either split or anisotropic over $K$. 

The surjectivity of word maps given by the word $w=[x,y]$ for different types of groups has a long history, see. \cite{EllersGordeevConj}, \cite{HuiLarShal}. For a split, or more generally quasi-split, group $\mathcal{G}^+(K)$ such surjectivity can be proved using a method called ''the Gauss decomposition with a given semi-simple part''. More precisely it can be proved for the group $\mathcal{G}^+(K)/Z(\mathcal{G}^+(K))$ if the field is assumed to be sufficiently large, see \cite{EllersGordeevConj}. Let us recall the main idea of this method, see \cite{CherEllersGordeev}, \cite{EllersGordeevGauss}, \cite{GordeevSumsofOrbits}, \cite{GordeevSaxl}, \cite{MoritaPlotkin}, \cite{Shengs}. Let $\mathcal{G}$ be a quasi-split group; $T$ be a fixed maximal quasi-split torus of $\mathcal{G}$; $B$ be a fixed Borel subgroup containing $T$; by $B^-$ we denote the opposite Borel subgroup. The method is based on the possibility for each noncentral conjugacy class of the group $\mathcal{G}^+(K)$ and any element $h\in T\cap \mathcal{G}^+(K)$ to find a representative $g$ of this class with Gauss decomposition $g=vhu$, where $u,v\in \mathcal{G}^+(K)$ are elements of unipotent radicals of subgroups $B$ and $B^-$ correspondently.

In paper \cite{KursovGL}, Kursov study the products of commutators in general linear group over a skew-field, see also \cite{KursovDiss}. The skew-field was assumed to be finite-dimensional over its centre and not coinciding with it. Later Kursov's result was rediscovered by Gordeev and Egorchenkova in paper \cite{CommInGLD}, see also \cite{EgorchDiss}. In this paper, the corresponding analog of the Gauss decomposition with a given semi-simple part was proved.

However, this analog does not allow to fully prescribe the diagonal factor $h$: the last element on the diagonal can be prescribed only up to multiplication by an element of the derived subgroup of the skew-field multiplicative group. Hence such a decomposition does not allow to prove that any noncentral element of the elementary group is a commutator. Nevertheless, it is proved in \cite{KursovGL} that if any element of the derived subgroup of the skew-field multiplicative group is a product of at most $c$ commutators, then any noncentral element of the elementary group is also a a product of at most $c$ commutators in general linear group. 

In fact, if there is no restriction on the commutator lengths in the derived subgroup of the skew-field multiplicative group, then it is impossible to express an arbitrary noncentral element of the elementary group as one commutator. 

It is well known that if an element $\tau\in D^*$ of the multiplicative group of the skew-field $D$ does not belong to the derived subgroup $[D^*,D^*]$, then the diagonal matrix $g=\diag(1,\ldots,1,\tau)$ does not belong to the elementary subgroup because it has a nontrivial Dieudonne determinant.

To a reader who is familiar with basic model theory, we left as an exercise to use the fact above and the compactness theorem to deduce the following statement: for any integers ${n\ge 2}$ and ${d\ge 1}$ there exists a positive integer $N$ such that for any skew-field  $D$ (not necessarily finite-dimensional over centre) and any $\tau\in D^*$, if $\tau$ is not a product of $N$ commutators in multiplicative group, then $\diag(1,\ldots,1,\tau)$ is not a product of $d$ commutators in general linear group of $n\times n$ matrices. In addition, if we have $\tau\in [D^*,D^*]$, then the matrix $\diag(1,\ldots,1,\tau)$ belong to the elementary subgroup. Therefore, if the derived subgroup of the multiplicative group contain an element of a large commutator length, then the elementary group also contain an element of a large commutator length.

In the present paper we prove an effective version of this statement. Namely we show that one can take $N=d(8n^2-13n+8)-2n^2+3n-1$ (Theorem~\ref{LengthOfTau}). Therefore, there is a lower estimate on the maximal commutator length in the elementary group (Corollary~\ref{lower}), that depend on maximal commutator length in derived subgroup of the multiplicative group. Hence there is a necessary condition for expressing an arbitrary noncentral element as a commutator (Corollary~\ref{necessary}).

However the main result of \cite{KursovGL} and \cite{CommInGLD} is not optimal. As we show, for skew-fields that are finite dimensional over its centre, if any element of the derived subgroup of the multiplicative group is a product of at most $c$ commutators, then any element of the elementary subgroup is a product of at most $\ceil{\tfrac{c}{n}}$ commutators in the general linear group of $n\times n$ matrices, and also a product of at most $\ceil{\tfrac{c}{n-2}}$ commutators in the elementary group, where by $\ceil{x}$ we denote the ceil function of $x$ (Theorem \ref{upper}). In particular, there is a sufficient condition for expressing an arbitrary noncentral element as a commutator (Corollaries \ref{sufficient1} and \ref{sufficient2}).
	
I am grateful to my teacher N. A. Vavilov for the thorough reading of this paper and valuable remarks. 

\section{Basic notation}
 
 \subsection{General linear group and elementary subgroup}
 
 Let $D$ be a skew-field, $D^*=D\sm\{0\}$ be its multiplicative group. We fix a positive integer $n\ge 2$. By $\GL(n,D)$ we denote the general linear group of degree $n$ over the skew-field $D$. By $\diag(\eps_1,\ldots,\eps_n)$ we denote the diagonal matrix with elements $\eps_1$,$\ldots$,$\eps_n$ on the diagonal. By $h_{i,j}(\eps)$, where $1\le i\ne j\le n$ and $\eps\in D^*$, we denote the diagonal matrix that has the element $\eps$ in $i$-th diagonal position, the element $\eps^{-1}$ in $j$-th diagonal position, and ones in all the remaining diagonal positions.
 
 We denote by $t_{i,j}(\xi)$, where $1\le i,j\le n$, $i\ne j$, $\xi\in D$, the elementary transvection, i.e. the matrix that has the element $\xi$ in the position $(i,j)$, and in all the other positions it coincides with the identity matrix.

 Let us also introduce notation for the root subgroups $X_{i,j}=t_{i,j}(D)\le\GL(n,D)$. 
 
 It is well known that the elementary transvections satisfy the following relations. 
 
$$
t_{i,j}(\xi)t_{i,j}(\zeta)=t_{i,j}(\xi+\zeta)\tc\eqno(1)
$$
$$
[t_{i,j}(\xi),t_{p,q}(\zeta)]=\begin{cases}
t_{i,q}(\xi\zeta)\tc\;&\text{ if } j=p\tc i\ne q\tc\\
e\tc\;&\text{ if }j\ne k\tc i\ne q\tc
\end{cases}\eqno(2)
$$
$$
t_{i,j}(\zeta)t_{j,i}(\xi)=h_{i,j}(\zeta)h_{i,j}(\zeta^{-1}+\xi)t_{j,i}((1+\xi\zeta)\xi)t_{i,j}((1+\zeta\xi)^{-1}\zeta)\tc\;\text{ if } \zeta, 1+\zeta\xi\in D^*\tc\eqno(3)
$$
$$
t_{i,j}(\xi)t_{j,i}(-\xi^{-1})t_{i,j}(\xi)=t_{j,i}(-\xi^{-1})t_{i,j}(\xi)t_{j,i}(-\xi^{-1})\tp\eqno(4)
$$

 The subgroup
  $$
  E(n,D)=\big\langle t_{i,j}(\xi)\colon 1\le i\ne j\le n\tc\; \xi\in D \big\rangle\tc
  $$ 
 generated by all the elementary transvections is called the {\it elementary group}. In case of a skew-field, the elementary group coincides with the kernel of the Dieudonne determinant $\map{\det}{\GL(n,D)}{D^*/[D^*,D^*]}$, see, for example, \cite{ArtinGeom}. In particular, we have $[\GL(n,D),\GL(n,D)]\le E(n,D)$. In fact, except for the case, where $n=2$ and $D=\F_2$, we have $E(n,D)=[\GL(n,D),\GL(n,D)]$.
 
 \subsection{Subgroups of triangular and diagonal matrices}
 
 By $U$,$V\le \GL(n,D)$ we denote the groups of upper and lower unitriangular matrices correspondently; i.e. triangular matrices with ones on the diagonal. Note that $U$,$V\le E(n,D)$.  
 
 Set 
 $$
 H=\{\diag(\eps_1,\ldots,\eps_n)\colon \eps_i\in D^*\}\le \GL(n,D)\tp
 $$
 
 Note that the Dieudonne determinant of the matrix $h=\diag(\eps_1,\ldots,\eps_n)$ is equal to the image of the product $\eps_1\ldots\eps_n$ in the quotient group $D^*/[D^*,D^*]$. Therefore, the element $h$ belongs to the elementary subgroup if and only if the product $\eps_1\ldots\eps_n$ belong to the derived subgroup of the group $D^*$.
 
 Set $B=HU\le\GL(n,D)$, it is the group of upper triangular matrices; and $B^-=HV\le\GL(n,D)$ is the group of lower triangular matrices.
 
 For $\kappa\in\N_0^{n-1}$, we denote by $H_{\le\kappa}$ the set of elements $h\in H$ that can be expressed as a product of elements $h_{i,i+1}(\eps)$, where $1\le i\le n-1$ and $\eps\in D^*$, in such a way that the number of factors with given $i$ is not greater than $\kappa_i$. The order of factors can be arbitrary.
 
 \subsection{Commutator lengths}
 
 Let $G$ be a group. We denote by $[G,G]_{\le k}$ the set of elements that are products of at most $k$ commutators. For an element $g\in [G,G]$, we denote by $l_G(g)$ its commutator length, i.e. the smallest integer $k$ such that $g\in [G,G]_{\le k}$. If the group in question is clear from the context, then we will write just $l(g)$.
 
 By $Z(G)$ we denote the centre of the group $G$.
 
Set 
 	$$
	c=\sup l_{D^*}(\xi)\tc\qquad \xi\in [D^*,D^*]\tp
	$$
	Except for the case $n=2$, $D=\F_2$, set
 	$$
 	C=\sup l_{\GL(n,D)}(g)\tc\qquad g\in E(n,D)\sm Z(E(n,D))\tp
 	$$
 	For $n\ge 3$, we also set
 	$$
 	C'=\sup l_{E(n,D)}(g)\tc\qquad g\in E(n,D)\sm Z(E(n,D))\tp
	$$

The aim of the present paper is to estimate the numbers $C$ and $C'$ for a given value of~$c$. 

Let us also note that if $g\in Z(E(n,D))$, then we can decompose it as a product $g=g_1g_2$, where $g_2$ is an arbitrary noncentral commutator in $E(n,D)$ (say, $g_2=[t_{12}(1),t_{21}(1)]$), and $g_1=gg_2^{-1}$ is a noncentral element. Therefore, the commutator length of the element $g$ in the group $\GL(n,D)$ (resp. $E(n,D)$) is not greater than $C+1$ (resp. $C'+1$).  
 
\section{A lower estimate}

In this section we prove the following theorem. 

{\thm\label{LengthOfTau} Let $D$ be an arbitrary skew-field. Let $\tau\in D^*$ be such that the matrix $\diag(1,1,\ldots,1,\tau)$ is a product of $d$ commutators in the group $\GL(n,D)$, where $d\ge1$. Then the element $\tau$ is a product of at most 
	$$
	d(8n^2-13n+8)-2n^2+3n-1
	$$
	 commutators in the group $D^*$. } 

From this theorem we obtain the following corollary.

{\cor\label{lower} In our notation, if $c\ge 1$, then the following inequality holds true:
$$
C\ge \frac{c+2n^2-3n+1}{8n^2-13n+8}\tp
$$
} 

The proof consists of several lemmas. 

Let us introduce the following notation. By $e_k\in\N_0^{n-1}$ we denote the row that has one in $k$-th position and zeroes in all the other positions. 

We denote by $U_k$ the subgroup of the group $U$ that consists of upper unitriangular matrices with zero in the position $(k,k+1)$. Using $(2)$, it is easy to see that $U$ is a semidirect product $U=U_k\leftthreetimes X_{k,k+1}$ and that $U_k^{t_{k+1,k}(\theta)}\le U$ for any $\theta\in D$. 

The group $V_k$ is defined similarly. Then $V=V_k\leftthreetimes X_{k+1,k}$.

{\lem\label{UVUt} Let $1\le k\le n-1$, let $\xi\in D$. Then the following inclusion holds true:
$$
UVUt_{k+1,k}(\xi)\le H_{\le 2e_k}UVU\tp
$$
\vspace*{-4ex}}
\begin{proof}
	
	Let $u_1,u_2\in U$ and $v\in V$. Let us prove that $u_1vu_2t_{k+1,k}(\xi)\in H_{\le 2e_k}UVU$. 
	
	The element $u_2$ can be written as $u_2=u_2't_{k,k+1}(\zeta)$, where $u_2'\in U_k$. Similarly, the element $v$ can be written as $v=v't_{k+1,k}(\eta)$, where $v'\in V_k$. 
	
	\smallskip
	
	{\bf Case 1}: $\zeta=0$. Using $(2)$, we obtain
	$$
	u_1vu_2t_{k+1,k}(\xi)=u_1vu_2't_{k+1,k}(\xi)=u_1vt_{k+1,k}(\xi)(u_2')^{t_{k+1,k}(\xi)}\in UVU\tp
	$$ 
	
	\smallskip
	
	{\bf Case 2}: $\zeta\in D^*$, but $\zeta\xi\ne -1$. Using $(3)$, we obtain
	\begin{multline*}
	u_1vu_2t_{k+1,k}(\xi)=u_1vu_2't_{k,k+1}(\zeta)t_{k+1,k}(\xi)=\\=u_1vu_2'h_{k,k+1}(\zeta)h_{k,k+1}(\zeta^{-1}+\xi)t_{k+1,k}((1+\xi\zeta)\xi)t_{k,k+1}((1+\zeta\xi)^{-1}\zeta)\tp
	\end{multline*}
	Further, using that the subgroup $H$ normalises the subgroups $U$,$V$ and $U'$ we can move the factors $h_{k,k+1}(\zeta_k)h_{k,k+1}(\zeta_k^{-1}+\xi)$ to the begining; after that we do the same transformation as in the first case and we obtain $u_1vu_2t_{k+1,k}(\xi)\in H_{\le 2e_k}UVU$.
	
	\smallskip
	
	{\bf Case 3}: $\eta\zeta\ne -1$. Acting as in the previous cases, we obtain that $u_1vt_{k,k+1}(\zeta)\in H_{\le 2e_k}UV$. Therefore,
	$$
	u_1vu_2=u_1vu_2't_{k,k+1}(\zeta)=u_1vt_{k,k+1}(\zeta)(u_2')^{t_{k,k+1}(\zeta)}\in H_{\le 2e_k}UV(u_2')^{t_{k,k+1}(\zeta)}\tc
	$$
	where $(u_2')^{t_{k,k+1}(\zeta)}\in U_k$. Further, applying calculation from the first case, we obtain $u_1vu_2t_{k+1,k}(\xi)\in H_{\le 2e_k}UVU$.
	
	\smallskip
	
	{\bf Case 4}: $\zeta\xi=-1$ and $\eta\zeta= -1$, i.e. $\eta=\xi$ and $\zeta=-\xi^{-1}$. In this case, using $(4)$, we obtain:
	\begin{multline*}
	u_1vu_2t_{k+1,k}(\xi)=u_1v't_{k+1,k}(\eta)u_2't_{k,k+1}(\zeta)t_{k+1,k}(\xi)=\\=u_1v't_{k+1,k}(\xi)t_{k,k+1}(-\xi^{-1})t_{k+1,k}(\xi)(u_2')^{t_{k,k+1}(-\xi^{-1})t_{k+1,k}(\xi)}=\\=u_1v't_{k,k+1}(-\xi^{-1})t_{k+1,k}(\xi)t_{k,k+1}(-\xi^{-1})(u_2')^{t_{k,k+1}(-\xi^{-1})t_{k+1,k}(\xi)}\tp
	\end{multline*}
	As in the firs case, we obtain that  $u_1v't_{k,k+1}(-\xi^{-1})\in UV$. Further $(u_2')^{t_{k,k+1}(-\xi^{-1})}\in U_k$; hence $(u_2')^{t_{k,k+1}(-\xi^{-1})t_{k+1,k}(\xi)}\in U$. Therefore, $u_1vu_2t_{k+1,k}(\xi)\in UVU$.
\end{proof}

Set $\lambda_1=2(n-1)$ and $\lambda_i=4(n-i)$ for $2\le i\le n-1$. 

{\lem\label{UVUV} The following inclusion holds true:
$$
UVUV\le H_{\le\lambda}UVU\tp
$$
\vspace*{-4ex}}
	\begin{proof}
		It is easy to see that any element of the group $V$ can be written as a product of elementary trancvections $t_{k+1,k}(\xi)$, among that there are at most $n-1$ transvectons with $k=1$, and at most $2(n-k_0)$ transvections with $k=k_0$ for $2\le k_0\le n-1$. Therefore the statement follows from Lemma~\ref{UVUt}.
	\end{proof}

{\lem\label{GLisHUVU} The following equality holds true:
$$
\GL(n,D)=HUVU\tp
$$
\vspace*{-4ex}}
\begin{proof}
	It is well known that $\GL(n,D)=H\cdot E(n,D)$; and the group $E(n,D)$ is generated by subgroups $U$ and $V$. It remains to apply Lemma~\ref{UVUV}. 
\end{proof}

Note that the last lemma is well known and holds true for semilocal rings, see, for example, \cite{BorevichParab}. Simillar result can also be found in \cite{SmolSuryVavGauss}.

Set $\mu=(6,3,\ldots,3)\in\N_0^{n-1}$ (in all positions, except the first one, it has three).

{\lem\label{CommInH} The following inclusion holds true:
$$
[H,H]_{\le 1}\le H_{\le\mu}
$$
\vspace*{-4ex}}
\begin{proof}
The statement follows from the following equality in $\GL(2,D)$:
$$
\diag([\xi,\zeta],1)=h_{1,2}(\xi)h_{1,2}(\zeta)h_{1,2}(\xi^{-1}\zeta^{-1})\tp
\vspace*{-3ex}
$$
\end{proof}

Set $\kappa^p=p\mu+(4p-1)\lambda$, where $p\in\N$. Therefore, $\kappa^1=\mu+3\lambda$ and $\kappa^p=\kappa^{p-1}+\kappa_1+\lambda$ for $p\ge 2$. 

{\lem\label{CommInGL} The following inclusion holds true: 
$$
[\GL(n,D),\GL(n,D)]_{\le p}\le H_{\le\kappa^p}UVU\tp
$$
\vspace*{-4ex}}
\begin{proof}
	The base of induction is $p=1$.  Using Lemma~\ref{GLisHUVU}, the fact that the subgroup $H$ normalises the subgroups $U$ and $V$, and also Lemmas~\ref{CommInH} and~\ref{UVUV}, we obtain
	$$
	[\GL(n,D),\GL(n,D)]_{\le 1}=[HUVU,HUVU]_{\le 1}\le [H,H]_{\le 1}UVUVUVUVU\le H_{\le\kappa^1}UVU\tp
	$$ 
	
	Now let us perform the induction step from $p-1$ to $p$. Using induction hypothesis, the base of induction, the fact that the subgroup $H$ normalises the subgroups $U$ and $V$, and also lemma~\ref{UVUV}, we obtain
	\begin{align*}
	[\GL(n,D),\GL(n,D)]_{\le p}=[\GL(n,D),\GL(n,D)]_{\le p-1}[\GL(n,D),\GL(n,D)]_{\le 1}\le\\ \le H_{\le\kappa^{p-1}}UVUH_{\le\kappa^1}UVU\le H_{\le\kappa^{p-1}+\kappa^1}UVUVU\le H_{\kappa^p}UVU\tp\qquad\qedhere
	\end{align*}
\end{proof}

{\lem\label{capH} For any $\kappa\in\N_0^{d-1}$ the following inclusion holds true:
$$
(H_{\le \kappa}UVU)\cap H=H_{\le \kappa}\tp
$$
\vspace*{-4ex}}
\begin{proof}
	\begin{multline*}
	(H_{\le \kappa}UVU)\cap H=(H_{\le \kappa}UVU)\cap B\cap H=(H_{\le \kappa}U(V\cap B)U)\cap H=\\=(H_{\le \kappa}U)\cap H=H_{\le \kappa}(U\cap H)=H_{\le \kappa}\tp\qedhere
	\end{multline*}
\end{proof}

Let us forget about general linear group for now, and prove two purely group-theoretic lemmas. 

Let $G$ be an arbitrary group. Recall that for any element $g\in [G,G]$ we denote by $l(g)$ its commutator length, i.e the smallest $k$ such that $g\in [G,G]_{\le k}$.

{\lem\label{permute} Let $a_1$,$\ldots$,$a_p$,$b_1$,$\ldots$,$b_q\in G$, where $q\ge 1$. Set $a=a_1^{-1}\ldots a_p^{-1}$, and $b=b_1^{-1}\ldots b_q^{-1}$. Let $w$ be a word with letters $a_1^{-1}$,$\ldots$,$a_p^{-1}$,$b_1$,$\ldots$,$b_q$, where each letter occur exactly once, and the letters $a_i^{-1}$ relatively to each other are arranged in indices ascending order. Assume that the element of the group $G$ represented by the word $w$ is equal to the identity element. Then if at least one of the elements $a$ and $b$ belongs to the derived subgroup of the group $G$, then the second element also belongs to the derived subgroup and the following inequality holds true: $l(b)\le l(a)+ q-1$. }

\begin{proof}
	
	We can reduce the problem to the case where the first letter in the word $w$ is $b_1$. In order to do this reduction we perform the corresponding cyclic permutation of letters in the word $w$, and after that we change the numeration of the elements $a_i$ (again by cyclic permutation), so that the condition remains true. The element $a$, therefore, became replaced by a conjugate, and its commutator length remains unchanged. 

	So assume that the letter $b_1$ is the first one. Consider the sequence of transformations of the word $w$ that has $q-1$ steps, where on $i$-th step we move the letter $b_{i+1}$ to the beginning without changing the relative order of the other letters. One such step can not increase the commutator length of the element represented by the word more than by one. After applying all the transformations we obtain the word that represents the element $b^{-1}a$. Therefore, $b^{-1}a\in [G,G]_{\le q-1}$, which implies the statement in question.   
	 
\end{proof}

{\lem\label{inverse} Let $a_1$,$\dots$,$a_k\in G$ be such that $a_1\ldots a_k=e$. Then the following inequality holds true $l(a_1^{-1}\ldots a_k^{-1})\le \max(0,k-2)$.}

\begin{proof}
	For $k=1$ the statement is obvious; so assume that $k\ge 2$.
	
	It is clear that $l(a_1^{-1}\ldots a_k^{-1})=l((a_1^{-1}\ldots a_k^{-1})^{-1})=l(a_k\ldots a_1)$. Further $l(a_k\ldots a_1)=l(a_1a_k\ldots a_2)$ because these elements are conjugate. Finally, the word $a_1a_k\ldots a_2$ can be obtained from the word $a_1\ldots a_k$ by $k-2$ transformations, where $i$-th transformation moves the letter $a_{k-i}$ to the end without changing the relative order of the other letters. One such step can not increase the commutator length of the element represented by the word more than by one. Therefore, $l(a_1a_k\ldots a_2)\le k-2$.
\end{proof}

Now let us return to the general linear group over a skew-field. Let $\kappa\in \N_0^{n-1}$. Set 
$$
s(\kappa)=\max(0,\kappa_{1}-2)+\sum_{i=2}^{n-1}\max(0,\kappa_i-1)\tp
$$

{\lem\label{skappa} Let $\kappa\in \N_0^{n-1}$, and let $\tau\in D^*$ be such that the matrix $g=\diag(1,\ldots,1,\tau)$ belongs to $H_{\le\kappa}$. Then $\tau\in [D^*,D^*]_{\le s(\kappa)}$.}
\begin{proof}
	If $\kappa_{n-1}=0$, then $\tau=1$, and the statement is trivial. Assume that $\kappa_{n-1}\ne 0$. Moreover, without loss of generality, we may assume that $\kappa_i\ne 0$ for all $i$. Indeed, otherwise consider the largest $k$ such that $\kappa_k=0$. It is easy to see that in this case we have
	$$
	\diag(1,\ldots,1,\tau)\in H_{\le (\kappa_{k+1},\ldots,\kappa_{n-1})}\le GL(n-k,D)\tc
	$$
	and we can prove the statement for this element.
	
	So let $\kappa_i\ne 0$ for all $i$. Let $g$ be the product of the elements $h_{i,i+1}(\xi_{i,j})$, $1\le j\le \kappa_i$ in some order. In addition, we assume that for a fixed $i$ these elements relatively to each other are arranged in index $j$ ascending order.
	
	Set $\eps_i=\xi_{i,1}^{-1}\ldots\xi_{i,\kappa_i}^{-1}\in D^*$. It follows from Lemma~\ref{inverse} and equality $g_{11}=1$ that $l(\eps_1)\le \max(0,\kappa_1-2)$. 
	
	For $2\le i\le n-1$ it follows from Lemma~\ref{permute} and equality $g_{i,i}=1$ that 
	$$
	l(\eps_i)\le l(\eps_{i-1})+\max(0,\kappa_i-1)\tp
	$$
	
	Therefore,
	 $$
	 l(\tau)=l(\eps_{n-1})=l(\eps_1)+\sum_{i=2}^{n-1}\left(l(\eps_i)-l(\eps_{i-1})\right)\le s(\kappa)\tp
	 \vspace*{-6ex}
	$$
\end{proof}

\bigskip

Now we can finish the proof of Theorem~\ref{LengthOfTau}. Let
 $$
 g=\diag(1,\ldots,1,\tau)\in [\GL(n,D),\GL(n,D)]_{\le d}\tp
 $$
 Using Lemmas~\ref{CommInGL}, \ref{capH} and~\ref{skappa}, we obtain that $\tau\in [D^*,D^*]_{\le s(\kappa^d)}$; and it is easy to see that
$$
s(\kappa^d)=d(8n^2-13n+8)-2n^2+3n-1\tp
$$

\section{An upper estimate}

In this section we prove th following theorem; its proof is the modification of the proof of the main theorem in \cite{CommInGLD}.

{\thm\label{upper} Let $D$ be a skew-field with the centre $K$. Let $1<\dim_KD<\infty$ (hence $K$ is an infinite field). Then in our notation the following inequality holds true:
	$$
	C\le \ceil*{\frac{c}{n}}\tc
	$$
	where by $\ceil{x}$ we denote the ceil function of $x$.
	 
	For $n\ge 3$ the following inequality also holds true:
	$$
	C'\le \ceil*{\frac{c}{n-2}}\tp
	$$
}

In the following two lemmas we assume the condition of Theorem~\ref{upper}.

 {\lem\label{conj} Let $c=d_1+\ldots+d_n$ be a partition of the number $c$, where $d_n\in \N_0$. Let $g\in E(n,D)\sm Z(E(n,D))$. Then there exists $\gamma\in E(n,D)$ such that $g^\gamma=vhu$, where $v\in V$, $u\in U$ and $h=\diag(\eps_1,\ldots,\eps_n)$, where $\eps_i\in [D^*,D^*]_{\le d_i}$ for all $1\le i\le n$.}
\begin{proof}
	Let $k$ be the smallest integer such that $d_k\ne 0$.
	
	The base of induction is for $k=n$, i.e. $d_n=c$. By Theorem 2.1. of \cite{CommInGLD} there exists an element $\gamma\in E(n,D)$ such that $g^\gamma=vhu$, where $v\in V$, $u\in U$ and $h=\diag(1,\ldots,1,\eps_n)$. Since $g\in E(n,D)$, it follows that $\eps_n\in [D^*,D^*]=[D^*,D^*]_{\le c}$. The base of induction is proved.
	
	Now let us perform the induction step from $k+1$ to $k$. By induction hypothesis there exists $\gamma_1\in E(n,D)$ such that $g^{\gamma_1}=v_1h_1u_1$, where $v_1\in V$, $u_1\in U$ and $h_1=\diag(1,\ldots,1,\tilde{\eps}_{k+1},\ldots,\tilde{\eps}_n)$, where $\tilde{\eps}_i\in [D^*,D^*]_{\le d_i}$ for $i\ge k+2$, and $\tilde{\eps}_{k+1}\in[D^*,D^*]_{\le d_k+d_{k+1}}$. 
	
	\smallskip
	
	{\bf Case 1}: $u_1=u_1't_{k,k+1}(\zeta)$, where $u_1'\in U_k$ and $\zeta\in D^*$. We express the element $\tilde{\eps}_{k+1}$ as a product $\tilde{\eps}_{k+1}=\eps_{k+1}\theta$, where $\eps_{k+1}\in [D^*,D^*]_{\le d_{k+1}}$ and $\theta\in[D^*,D^*]_{\le d_k}$.
	
	 Since $\zeta\in D^*$, it follows that $\theta=1+\xi\zeta$ for some $\xi\in D$. Set $\gamma=\gamma_1t_{k+1,k}(\xi)$. Then, applying $(3)$, it is easy to see that $g^\gamma=vhu$, where $v\in V$, $u\in U$ and 
	 $$
	 h=\diag(1,\ldots,1, 1+\zeta\xi,\eps_{k+1},\tilde{\eps}_{k+2},\ldots,\tilde{\eps}_n)\tp
	 $$
	 It remains to notice that the element $1+\zeta\xi$ is conjugate to the element $1+\xi\zeta=\theta$; hence it has the same commutator length. 

	\smallskip
		 
	 {\bf Case 2}: $v_1=v_1't_{k+1,k}(\zeta)$, where $v_1'\in V_k$ and $\zeta\in D^*$. The proof for this case is similar to the previous one.
	 
	\smallskip
		 
	 {\bf Case 3}: $u_1\in U_k$ and $v_1\in V_k$. It is clear that if $\tilde{\eps}_{k+1}=1$, then one can take $\gamma=\gamma_1$. Now assume that $\tilde{\eps}_{k+1}\ne 1$. Set $\gamma_2=\gamma_1t_{k,k+1}(1)$. Then
	 $$
	 g^{\gamma_2}=t_{k,k+1}(-1)v_1h_1u_1t_{k,k+1}(1)=v_1^{t_{k,k+1}(1)}h_1[h_1^{-1},t_{k,k+1}(-1)]u_1^{t_{k,k+1}(1)}\tc
	 $$
	 where $v_1^{t_{k,k+1}(1)}\in V$, $u_1^{t_{k,k+1}(1)}\in U_k$; and since $\tilde{\eps}_{k+1}\ne 1$, it follows that $[h_1^{-1},t_{k,k+1}(-1)]=t_{k,k+1}(\zeta)$, where $\zeta\in D^*$. Therefore, the problem is reduced to the first case.
\end{proof}

{\lem\label{OneComm} Let $g=vhu\in \GL(n,D)$, where $v\in V$, $u\in U$ and $h=\diag(\eps_1,\ldots,\eps_n)$, where $\eps_i\in [D^*,D^*]_{\le 1}$. Then $g\in [\GL(n,D),\GL(n,D)]_{\le 1}$.
	
	 If, in addition, the elements $\eps_1$ and $\eps_2$ are equal to one, then $g\in [E(n,D),E(n,D)]_{\le 1}$.}
\begin{proof}
	Let $\eps_i=[a_i,b_i]$, where $a_i$,$b_i\in D^*$. If $\eps_1=\eps_2=1$, then we may assume that $a_2$, $b_1\in K$, $a_1=(a_2\ldots a_n)^{-1}$, and $b_2=(b_1b_3\ldots b_n)^{-1}$.
	
	Since the field $K$ is infinite, it follows that replacing for all $i\ge 2$ the elements $a_i$ by $t_ia_i$ and also replacing $a_1$ by $(t_2\ldots t_n)^{-1}a_1$ for suitable $t_i\in K^*$ we may assume that the reduced norms of the elements $a_i$ are pairwise distinct. The elements $\eps_i$ remains the same after such a replacing. The condition $a_1=(a_2\ldots a_n)^{-1}$ is preserved if it were assumed.
	 
	Set 
	 \begin{gather*}
	 h_1=\diag(a_1,\ldots,a_n)\tc\qquad
	 \tau=\diag(b_1,\ldots,b_n)\tc\\
	 h_2=\tau h_1^{-1}\tau^{-1}=\diag(b_1a_1^{-1}b^{-1},\ldots,b_na_n^{-1}b_n^{-1})\tp
	 \end{gather*}
	Therefore $h=h_1h_2$. Since the reduced norms of the elements $a_i$, and hence of the elements $b_ia_ib_i^{-1}$, are pairwise distinct, it follows that the elements $h_1$ and $h_2$ are $D$-regular (see Definition 3.3 and Proposition 3.4 in \cite{CommInGLD}). By Proposition 4.1 in \cite{CommInGLD} there exist $u'\in U$ and $v'\in V$ such that $v=[v',h_1]$ and $u=[h_2^{-1},u']$. 
	  
	Similarly to the proof of the main theorem in \cite{CommInGLD}, we obtain that
	\begin{multline*}
	g=vhu=v'h_1v'^{-1}u'h_2u'^{-1}=v'h_1v'^{-1}u'\tau h_1^{-1}\tau^{-1}u'^{-1}=\\=[v'h_1v'^{-1},u'\tau v^{-1}]\in[\GL(n,D),\GL(n,D)]_{\le 1}\tp
	\end{multline*}
	  Moreover, if $\eps_1=\eps_2=1$, then according to our assumptions $h_1$,$\tau\in E(n,D)$; hence $g\in[E(n,D),E(n,D)]_{\le 1}$.
\end{proof}

Let us finish the proof of Theorem \ref{upper}. Let $g\in E(n,D)\sm Z(E(n,D))$. We prove that $g\in [\GL(n,D),\GL(n,D)]_{\ceil{\frac{c}{n}}}$.

By Lemma \ref{conj} we may assume that $g=vhu$, where $v\in V$, $u\in U$ and $h=\diag(\eps_1,\ldots,\eps_n)$, where $\eps_i\in [D^*,D^*]_{\le \ceil{\frac{c}{n}}}$ for all $1\le i\le n$. For any $i$ let $\eps_i=\eps_i'\eps_i''$, where $\eps_i'\in [D^*,D^*]_{\le \ceil{\frac{c}{n}}-1}$ and $\eps_i''\in [D^*,D^*]_{\le 1}$. Set $h'=\diag(\eps_1',\ldots',\eps_n')$ and $h''=\diag(\eps_1'',\ldots,\eps_n'')$. We also set $\tilde{v}=v^{h'}$. 

Therefore, we have $g=vhu=h'\tilde{v}h''u$. By Lemma \ref{OneComm} we obtain that $\tilde{v}h''u\in[\GL(n,D),\GL(n,D)]_{\le 1}$. It is also clear that $h'\in [\GL(n,D),\GL(n,D)]_{\le \ceil{\frac{c}{n}}-1}$. Therefore, $g\in [\GL(n,D),\GL(n,D)]_{\ceil{\frac{c}{n}}}$.

Let $n\ge 3$. Let us prove that $g\in [E(n,D),E(n,D)]_{\ceil{\frac{c}{n-2}}}$.

By Lemma \ref{conj} we may assume that $g=vhu$, where $v\in V$, $u\in U$ and $h=\diag(1,1,\eps_3,\ldots,\eps_n)$, where $\eps_i\in [D^*,D^*]_{\le \ceil{\frac{c}{n-2}}}$ for all $3\le i\le n$. Similarly to the previous calculation, for any $i$ let $\eps_i=\eps_i'\eps_i''$, where $\eps_i'\in [D^*,D^*]_{\le \ceil{\frac{c}{n-2}}-1}$ and $\eps_i''\in [D^*,D^*]_{\le 1}$. Set $h'=\diag(1,1\eps_3',\ldots,\eps_n')$ and $h''=\diag(1,1,\eps_3'',\ldots,\eps_n'')$. We also set $\tilde{v}=v^{h'}$. 

As before, $g=h'\tilde{v}h''u$. By Lemma \ref{OneComm} we have $\tilde{v}h''u\in[E(n,D),E(n,D)]_{\le 1}$. It remains to show that $h'\in [E(n,D),E(n,D)]_{\le \ceil{\frac{c}{n-2}}-1}$.

For simplicity we set $d=\ceil{\frac{c}{n-2}}-1$. Let $\eps_i'=[a_{i,1},b_{i,1}]\ldots [a_{i,d},b_{i,d}]$. Set 
\begin{gather*}
h_{a,k}=\diag((a_{3,k}\ldots a_{n,k})^{-1},1,a_{3,k},\ldots,a_{n,k})\tc\\
h_{b,k}=\diag(1,(b_{3,k}\ldots b_{n,k})^{-1},b_{3,k},\ldots,b_{n,k})\tp
\end{gather*}

Therefore, $h'=[h_{a,1},h_{b,1}]\ldots [h_{a,d},h_{b,d}]\in [E(n,D),E(n,D)]_{\le d}$.

\section{Expressing matrices as one commutator}

Consider the following corollaries of our results concerning the question on expressing any noncentral element of the group $E(n,D)$ as one commutator.

Firstly, from the Corollary \ref{lower} we obtain the following corollary.

{\cor\label{necessary} Let $D$ be an arbitrary skew-field. Then for any noncentral element of the group $E(n,D)$ to be one commutator in $\GL(n,D)$, it is necessary that any element of the derived subgroup of the group $D^*$ were a product of at most $6n^2-10n+7$ commutators.}

Secondly, from Theorem \ref{upper} we obtain the following corollaries.

{\cor\label{sufficient1} Let $D$ be a skew-field and $K$ be its centre. Let $1<\dim_KD<\infty$. Then for any noncentral element of the group $E(n,D)$ to be one commutator in $\GL(n,D)$, it is sufficient that any element of the derived subgroup of the group $D^*$ were a product of at most $n$ commutators.}

{\cor\label{sufficient2} Let $D$ be a skew-field and $K$ be its centre. Let $1<\dim_KD<\infty$, and let $n\ge 3$. Then for any noncentral element of the group $E(n,D)$ to be one commutator in $E(n,D)$, it is sufficient that any element of the derived subgroup of the group $D^*$ were a product of at most $n-2$ commutators.}

\medskip

We denote by $\GL(D)$ the {\it stable general linear group}, i.e. $\GL(D)=\varinjlim \GL(n,D)$, where the limit is taken with respect to the inclusions.
$$
\GL(n,D)\to\GL(n+1,D)\tc\qquad
g\mapsto\begin{pmatrix}
g & 0\\
0 & 1
\end{pmatrix}\tp
$$
The subgroup $E(D)\le\GL(D)$ is defined as generated by elementary transvections. In other words, $E(D)=\varinjlim E(n,D)$.

In case where $c<\infty$, the following corollary follows directly from the previous one. However, using the proof of Theorem \ref{upper}, we can see that the condition $c<\infty$ can be lifted.

{\cor Let $D$ be a skew-field and $K$ be its centre. Let $1<\dim_KD<\infty$. Then any element of the group $E(D)$ is one commutator in $E(D)$.}  

\begin{proof}
	Let $g\in E(D)$. Since the centre of the group $\GL(D)$ is trivial, it follows that there exists $n>2$ such that $g\in E(n,D)\sm Z(E(n,D))$. 
	
	By Theorem 2.1. in \cite{CommInGLD}, which we refer in the proof of Lemma \ref{conj}, we may assume that $g=vhu$, where $v\in V$, $u\in U$ and $h=\diag(1,\ldots,1,\tau)$. Since $g\in E(n,D)$, it follows that $\tau\in [D^*,D^*]$; hence $\tau\in[D^*,D^*]_{\le d}$ for some $d$.
	
	Increasing the number $n$, if necessary, and conjugating $g$ by a permutation matrix in such a way that the previous assumption remains true, we may assume that $n\ge d+2$.
	
	In this case, it follows from the proof of Lemma \ref{conj} that the element $g$ is conjugate to the element $g_1=v_1h_1u_1$, where $v_1\in V$, $u_1\in U$ and $h_1=\diag(1,1,\eps_3\ldots,\eps_n)$, where $\eps_i\in [D^*,D^*]_{\le 1}$.
	
	Then by Lemma \ref{OneComm}, the element $g_1$, and hence the element $g$, is one commutator in $E(n,D)$.
\end{proof}

\providecommand{\bysame}{\leavevmode\hbox to3em{\hrulefill}\thinspace}
\providecommand{\MR}{\relax\ifhmode\unskip\space\fi MR }
\providecommand{\MRhref}[2]{%
	\href{http://www.ams.org/mathscinet-getitem?mr=#1}{#2}
}
\providecommand{\href}[2]{#2}

\end{document}